# Journal of Mathematics

# A Novel Necessary and Sufficient Condition for the Positivity of a Binary Quartic Form

Yang Guo

Department of Mathematics, Northeastern University, Shenyang, Liaoning 110819, China

Correspondence should be addressed to Yang Guo; guoyang@mail.neu.edu.cn

## Abstract

In this paper, by considering the common points of two conics instead of the roots of the binary quartic form, we propose a novel necessary and sufficient condition for the positivity of a binary quartic form using the theory of the pencil of conics. First, we show the degenerate members of the pencil of conics according to the distinct natures of the common points of two base conics. Then, the inequalities about the parameters of the degenerate members are obtained according to the properties of the degenerate conics. Last, from the inequalities we derive a novel criterion for determining the positivity of a binary quartic form without the discriminant.

## 1. Introduction

In $\mathbb{R}[x, y]$, let the polynomial

$$f(x,y) = x^4 + a_3 x^3 y + a_2 x^2 y^2 + a_1 x y^3 + a_0 y^4 \tag{1}$$

be a monic binary quartic form on real number field $\mathbb{R}$. We say that $f(x, y)$ is positive definite (positive semi-definite) if for any real numbers $x, y$ (not both 0), we have $f(x, y) > 0$ ($f(x, y) \geq 0$).

The problem of determining the global positivity or nonnegativity of forms has many applications. For instance, in the study of two-dimensional nonlinear systems and other mathematical and physical fields, we often need to deal with the positive definite or negative definite problems of a binary quartic form [1, 2]. So far, about the positive definite problem of a binary quartic form, there has existed some literatures. For example, Chin [3] used the canonical forms of a real quartic equation to propose a necessary and sufficient condition for no real roots. Gadenz and Li [4] published a systematic method for determining the positive definiteness of the binary quartic form. Their method adopted the testing of the permanence and variations of sign in the sequence of principal minors of the corresponding Hankel matrix. Jury and Mansour [5] used Ferrari's solution to obtain condition for positivity for a quartic equation. In the same year, Fuller [6] also gave similar criteria in terms of bigradient determinants. Ku [7] gave a criterion for the positive definiteness of the general quartic form from the known solution. However, his result was not complete. Wang and Qi [8] gave a complete and improved result.

In the above literatures, one of the most common results obtained by different methods can be summarized as follows. Let

$$V(x, y) = c_0 x^4 + 4c_1 x^3 y + 6c_2 x^2 y^2 + 4c_3 x y^3 + c_4 y^4 \qquad (2)$$

be a binary quartic form on real number field $\mathbb{R}$. In terms of the coefficients of (2), the following quantities can be defined [7]:

$$G = c_0^2 c_3 - 3c_0 c_1 c_2 + 2c_1^3;$$

$$H = c_0 c_2 - c_1^2;$$

$$I = c_0 c_4 - 4c_1 c_3 + 3c_2^2;$$

$$J = \det\left(\begin{bmatrix} c_0 & c_1 & c_2 \\ c_1 & c_2 & c_3 \\ c_2 & c_3 & c_4 \end{bmatrix}\right);$$

and

$$\Delta = I^3 - 27J^2.$$

Suppose $c_0 > 0$. Then the quartic form $V(x, y)$ is positive definite if and only if one of the following three conditions holds [3, 5, 6, 8]:

1) $\Delta = 0, G = 0, 12H^2 - c_0^2 I = 0, H > 0$;
2) $\Delta > 0, H \geq 0$;
3) $\Delta > 0, H < 0, 12H^2 - c_0^2 I < 0$.

Compared with the criterion in [4], this criterion based on the discriminant of a quartic polynomial is superior in that the conditions contained in this criterion are stated in the form of a set of inequalities which are explicit functions of the coefficients of the quartic form. Besides for this common criterion, Hasan et al. [9] presented an alternative criterion. They made use of some of the results of the quadratic programming theory to show that testing for positive semi-definiteness of the quartic form (1) was reduced to a test whether there is a real number $\lambda$ such that the parametric matrix

$$\mathbf{M}(\lambda) = \begin{bmatrix} 1 & \dfrac{a_3}{2} & -\dfrac{\lambda}{2} \\ \dfrac{a_3}{2} & a_2 + \lambda & \dfrac{a_1}{2} \\ -\dfrac{\lambda}{2} & \dfrac{a_1}{2} & a_0 \end{bmatrix}$$

is positive semi-definite. Essentially, Hasan et al.'s method needs to solve the inequality equations coming from the constraints that the parametric matrix $\mathbf{M}(\lambda)$ is positive semi-definite (positive definite).

Now, the natural question to ask is whether there exists a concrete real number $\lambda$ obtained by the coefficients of the form (1) such that the form (1) is positive semi-definite (positive definite) if and only if the matrix $\mathbf{M}(\lambda)$ is a positive semi-definite (positive definite) matrix. If yes, the criterion in [9] will be simplified. In this paper, we will see that the answer to this question is yes. We use the relationship between the common points of two conics and the degenerate members of the pencil of conics to address this question and further derive a novel criterion for determining the positivity of a binary quartic form without the discriminant. We show that given a monic binary quartic form (1), the following quantities can be defined:

$$b_0 = \frac{1}{4}(-a_1^2 + a_1 a_2 a_3 - a_0 a_3^2); \qquad b_1 = \frac{1}{4}(4a_0 - a_2^2 - a_1 a_3); \qquad b_2 = \frac{a_2}{2};$$

$$\lambda_0 = \frac{4b_2 + 2\sqrt{3b_1 + 4b_2^2}}{3}.$$

Then the form (1) is positive semi-definite (positive definite) if and only if the matrix

$$\begin{bmatrix} 1 & \frac{a_3}{2} & \frac{a_2 - \lambda_0}{2} \\ \frac{a_3}{2} & \lambda_0 & \frac{a_1}{2} \\ \frac{a_2 - \lambda_0}{2} & \frac{a_1}{2} & a_0 \end{bmatrix}$$

is a positive semi-definite (positive definite) matrix (Theorem 1).

More concretely, the form (1) is positive semi-definite if and only if

$$\lambda_0 \geq \frac{a_3^2}{4}, \qquad -\frac{1}{4}\lambda_0^3 + b_2 \lambda_0^2 + b_1 \lambda_0 + b_0 \geq 0.$$

In particular, the form (1) is positive definite if and only if (Theorem 2)

$$\lambda_0 > \frac{a_3^2}{4}, \qquad -\frac{1}{4}\lambda_0^3 + b_2 \lambda_0^2 + b_1 \lambda_0 + b_0 > 0.$$

## 2. Notations and Preliminaries

Throughout this paper, let $\mathbf{A}_1$ and $\mathbf{A}_2$ denote respectively two real symmetric matrices as follows:

$$\mathbf{A}_1 = \begin{bmatrix} 1 & \frac{a_3}{2} & \frac{a_2}{2} \\ \frac{a_3}{2} & 0 & \frac{a_1}{2} \\ \frac{a_2}{2} & \frac{a_1}{2} & a_0 \end{bmatrix}, \qquad \mathbf{A}_2 = \begin{bmatrix} 0 & 0 & -\frac{1}{2} \\ 0 & 1 & 0 \\ -\frac{1}{2} & 0 & 0 \end{bmatrix},$$

and let the matrix $\mathbf{M}_\lambda$ be

$$\mathbf{M}_\lambda = \mathbf{A}_1 + \lambda \mathbf{A}_2 = \begin{bmatrix} 1 & \dfrac{a_3}{2} & \dfrac{a_2 - \lambda}{2} \\ \dfrac{a_3}{2} & \lambda & \dfrac{a_1}{2} \\ \dfrac{a_2 - \lambda}{2} & \dfrac{a_1}{2} & a_0 \end{bmatrix},$$

where $\lambda$ is a variable real parameter. Since every $3 \times 3$ real symmetric matrix denotes a conic in the complex projective plane, $\mathbf{M}_\lambda$ is a pencil of conics determined by the base conics $\mathbf{A}_1$, $\mathbf{A}_2$.

Obviously, using the matrix $\mathbf{M}_\lambda$, the binary quartic form (1) can be represented as follows:

$$f(x, y) = [x^2, xy, y^2] \mathbf{M}_\lambda \begin{bmatrix} x^2 \\ xy \\ y^2 \end{bmatrix}.$$

On the other hand, for arbitrary $x$, $y$ not all zero, $[x^2, xy, y^2]^T$ is the projective coordinate of the point of the conic $\mathbf{A}_2$, so we can use the relative position of the conics $\mathbf{A}_1$, $\mathbf{A}_2$ to address the determination of positivity of the binary quartic form (1). We assert that the binary quartic form (1) is positive semi-definite if and only if there is no real simple (or transversal) intersection point between the conic $\mathbf{A}_1$ and the conic $\mathbf{A}_2$. In particular, the form (1) is positive definite if and only if there is no real common point between $\mathbf{A}_1$ and $\mathbf{A}_2$.

It is well known that the conics $\mathbf{A}_1$, $\mathbf{A}_2$ have four common points (which need not all be distinct). Further the relative position of the conics $\mathbf{A}_1$, $\mathbf{A}_2$ has nine cases according to the nature (generic or non-generic, real or complex) of their common points [10, 11]. They are presented in Table 1.

Table 1 Nine cases of the intersection between the conic $\mathbf{A}_1$ and the conic $\mathbf{A}_2$

| Case | Nature of the common point |
|------|---------------------------|
| 1 | four real simple points |
| 2 | two pairs of complex conjugate simple points |
| 3 | two real simple points + a pair of complex conjugate simple points |
| 4 | two real simple points + a real double point (simple-contact) |
| 5 | a pair of complex conjugate simple points + a real double point (simple-contact) |
| 6 | two real double points (double contact) |
| 7 | a pair of complex conjugate double points (double contact) |
| 8 | a real simple point + a triple point (three-point contact) |
| 9 | a quadruple point (four-point contact) |

Now, we discuss the degenerate members of the pencil $\mathbf{M}_\lambda$ since they are closely related to the relative position of the conics $\mathbf{A}_1$, $\mathbf{A}_2$. The pencil $\mathbf{M}_\lambda$ contains three degenerate members (which need not all be distinct) [12]. Their parameters are the roots of the cubic equation $\det(\mathbf{M}_\lambda) = 0$. The degenerate conics are the line-pairs through four common points of $\mathbf{A}_1$ and $\mathbf{A}_2$. So, according to the nature of the common point in Table 1, the degenerate member in the pencil $\mathbf{M}_\lambda$ may concretely consist of one of as follows:

(1) a real line-pair (two real lines);
(2) a complex line-pair (two complex lines);
(3) a complex conjugate line-pair (two complex conjugate lines);
(4) a real repeated line.

Table 2 Degenerate conics in the pencil $\mathbf{M}_\lambda$ and properties of three roots of $\det(\mathbf{M}_\lambda) = 0$

| Case [a] | Degenerate conic | Property of $\lambda_1, \lambda_2, \lambda_3$ |
|---|---|---|
| 1 | three real line-pairs | $\lambda_1 < \lambda_2 < \lambda_3 \leq \frac{a_3^2}{4}$ |
| 2 | one real line-pair + two complex conjugate line-pairs | $\lambda_1 \leq \frac{a_3^2}{4} < \lambda_2 < \lambda_3$ |
|   |   | or $\lambda_1 < \frac{a_3^2}{4} \leq \lambda_2 < \lambda_3$ |
| 3 | one real line-pair + two complex line-pairs | $\lambda_1 \leq \frac{a_3^2}{4}$, $\lambda_2 = \bar{\lambda}_3$ [b] |
| 4 | two real line-pairs | $\lambda_1 = \lambda_2 < \lambda_3 \leq \frac{a_3^2}{4}$ |
|   |   | or $\lambda_1 < \lambda_2 = \lambda_3 < \frac{a_3^2}{4}$ |
| 5 | one real line-pair + one complex conjugate line-pair | $\lambda_1 \leq \frac{a_3^2}{4} < \lambda_2 = \lambda_3$ |
| 6 | one real repeated line + one real line-pair | $\lambda_1 < \frac{a_3^2}{4} = \lambda_2 = \lambda_3$ |
| 7 | one real repeated line + one complex conjugate line-pair | $\lambda_1 = \lambda_2 = \frac{a_3^2}{4} < \lambda_3$ |
| 8 | one real line-pair | $\lambda_1 = \lambda_2 = \lambda_3 < \frac{a_3^2}{4}$ |
| 9 | one real repeated line | $\lambda_1 = \lambda_2 = \lambda_3 = \frac{a_3^2}{4}$ |

[a] Corresponding to Table 1

[b] $\lambda_2, \lambda_3$ are complex conjugate

The middle column of Table 2 gives the concrete type of every degenerate member in nine cases. The different type of degenerate member $\mathbf{M}_{\lambda_i}$ can provide important information about its parameter $\lambda_i$. More formally, we have the following lemmas.

**Lemma 1.** *If the matrix $\mathbf{M}_{\lambda_i}$ is a degenerate conic consisting of a real line-pair, then the rank of $\mathbf{M}_{\lambda_i}$ is 2 and $\lambda_i \leq \frac{a_3^2}{4}$.*

**Lemma 2.** *If the matrix $\mathbf{M}_{\lambda_i}$ is a degenerate conic consisting of a complex conjugate line-pair, then the rank of $\mathbf{M}_{\lambda_i}$ is 2 and $\lambda_i \geq \frac{a_3^2}{4}$.*

**Lemma 3.** *The matrix $\mathbf{M}_{\lambda_i}$ is a degenerate conic consisting of a real repeated line if and only if $\lambda_i = \frac{a_3^2}{4}$ is a root of multiplicity $\geq 2$ of $\det(\mathbf{M}_\lambda) = 0$.*

The first two lemmas can be proved by using Sylvester's criterion [13, Theorem 7.2.5] and the facts that the matrix $\mathbf{M}_{\lambda_i}$ is a degenerate conic consisting of a real line-pair if and only if $\mathbf{M}_{\lambda_i}$ is indefinite with the rank 2 and the matrix $\mathbf{M}_{\lambda_i}$ is a degenerate conic consisting of a complex conjugate line-pair if and only if $\mathbf{M}_{\lambda_i}$ is positive semi-definite with the rank 2. Lemma 3 can be proved by using direct computation and the fact that the matrix $\mathbf{M}_{\lambda_i}$ is a degenerate conic consisting of a real repeated line if and only if $\mathbf{M}_{\lambda_i}$ is positive semi-definite with the rank 1.

From the above three lemmas, it is easy to obtain the relation between $a_3^2/4$ and the parameters $\lambda_1, \lambda_2, \lambda_3$ of the three degenerate members of the pencil $\mathbf{M}_\lambda$ in every case that is given in the middle column of Table 2. The results are presented in the last column of Table 2. It is worth noting that this is a one-to-one correspondence between the middle column of Table 2 and the last column of Table 2.

We denote by $g(\lambda)$ the determinant of the parameter matrix $\mathbf{M}_\lambda$, namely

$$g(\lambda) = \det(\mathbf{M}_\lambda) = \det(\mathbf{A}_1 + \lambda \mathbf{A}_2) = -\frac{1}{4}\lambda^3 + b_2\lambda^2 + b_1\lambda + b_0,$$

where $b_0 = (-a_1^2 + a_1 a_2 a_3 - a_0 a_3^2)/4$, $b_1 = (4a_0 - a_2^2 - a_1 a_3)/4$, $b_2 = a_2/2$. Using the resultant of $g(\lambda)$ and its formal derivative $g'(\lambda)$, we may obtain the discriminant of $g(\lambda)$ as follows [14]:

$$D\big(g(\lambda)\big) = \frac{16(3b_1 + 4b_2^2)^3 - (27b_0 + 36b_1 b_2 + 32b_2^3)^2}{432}.$$

Since $D\big(g(\lambda)\big) < 0$ if and only if $g(\lambda) = 0$ has a pair of complex conjugate roots, we assert that $3b_1 + 4b_2^2 \geq 0$ for the other eight cases except for Case 3 in Table 2.

Now let

$$\lambda_0 = \frac{4b_2 + 2\sqrt{3b_1 + 4b_2^2}}{3}.$$

Obviously, when $3b_1 + 4b_2^2 \geq 0$, $\lambda_0$ is the largest of the two real roots (which need not be distinct) of $g'(\lambda) = 0$, otherwise $\lambda_0$ is a complex number.

Further, according to $g(-\infty) > 0$, $g(+\infty) < 0$ and the inequalities about three roots $\lambda_1, \lambda_2, \lambda_3$ in Table 2, we may get the properties of $\lambda_0$ and $g(\lambda_0)$ in the nine cases. For example, for $\lambda_1 \leq a_3^2/4$, $\lambda_2 = \bar{\lambda}_3$ in Case 3, i.e., $D\big(g(\lambda)\big) < 0$, then one of the following two configurations happens:

(1) if $3b_1 + 4b_2^2 < 0$, then $\lambda_0$ is a non-real number;
(2) if $3b_1 + 4b_2^2 \geq 0$, then $\lambda_0 < a_3^2/4$ or $\lambda_0 \geq a_3^2/4$, $g(\lambda_0) < 0$.

For the other eight cases in Table 2, $\lambda_0$ must be real and furthermore the relation between $\lambda_0$ and $a_3^2/4$ and the relation between $g(\lambda_0)$ and 0 can be obtained easily. The results are presented in the last column of Table 3.

Table 3 Properties of $\lambda_0$ and $g(\lambda_0)$

| Case [a] | Property of $\lambda_1, \lambda_2, \lambda_3$ | Property of $\lambda_0$ and $g(\lambda_0)$ |
|---|---|---|
| 1 | $\lambda_1 < \lambda_2 < \lambda_3 \leq \frac{a_3^2}{4}$ | $\lambda_0 < \frac{a_3^2}{4}$, $g(\lambda_0) > 0$ |
| 2 | $\lambda_1 \leq \frac{a_3^2}{4} < \lambda_2 < \lambda_3$ or $\lambda_1 < \frac{a_3^2}{4} \leq \lambda_2 < \lambda_3$ | $\lambda_0 > \frac{a_3^2}{4}$, $g(\lambda_0) > 0$ |
| 3 | $\lambda_1 \leq \frac{a_3^2}{4}$, $\lambda_2 = \bar{\lambda}_3$ | $\lambda_0$ is non-real or $\lambda_0 < \frac{a_3^2}{4}$ or $\lambda_0 \geq \frac{a_3^2}{4}$, $g(\lambda_0) < 0$ |
| 4 | $\lambda_1 = \lambda_2 < \lambda_3 \leq \frac{a_3^2}{4}$ or $\lambda_1 < \lambda_2 = \lambda_3 < \frac{a_3^2}{4}$ | $\lambda_0 < \frac{a_3^2}{4}$, $g(\lambda_0) \geq 0$ |
| 5 | $\lambda_1 \leq \frac{a_3^2}{4} < \lambda_2 = \lambda_3$ | $\lambda_0 > \frac{a_3^2}{4}$, $g(\lambda_0) = 0$ |
| 6 | $\lambda_1 < \frac{a_3^2}{4} = \lambda_2 = \lambda_3$ | $\lambda_0 = \frac{a_3^2}{4}$, $g(\lambda_0) = 0$ |
| 7 | $\lambda_1 = \lambda_2 = \frac{a_3^2}{4} < \lambda_3$ | $\lambda_0 > \frac{a_3^2}{4}$, $g(\lambda_0) > 0$ |
| 8 | $\lambda_1 = \lambda_2 = \lambda_3 < \frac{a_3^2}{4}$ | $\lambda_0 < \frac{a_3^2}{4}$, $g(\lambda_0) = 0$ |
| 9 | $\lambda_1 = \lambda_2 = \lambda_3 = \frac{a_3^2}{4}$ | $\lambda_0 = \frac{a_3^2}{4}$, $g(\lambda_0) = 0$ |

[a] Corresponding to Table 1 and Table 2

## 3. The Main Results

From the foregoing discussion, we know that the binary quartic form (1) is positive semi-definite if and only if there is no real simple intersection point between the conic $\mathbf{A}_1$ and the conic $\mathbf{A}_2$, i.e., Case 2, Case 5, Case 6, Case 7 and Case 9 in Table 1. In particular, the form (1) is positive definite if and only if there is no real common point between the conic $\mathbf{A}_1$ and the conic $\mathbf{A}_2$, i.e., Case 2, Case 7. Consequently, we can use the properties of $\lambda_0$ and $g(\lambda_0)$ in Table 3 to obtain a simplified necessary and sufficient condition compared with [9] for determining the positivity of the form (1).

**Theorem 1.** *Given a monic binary quartic form*

$$f(x,y) = x^4 + a_3 x^3 y + a_2 x^2 y^2 + a_1 x y^3 + a_0 y^4. \tag{3}$$

*Let*

$$b_1 = \frac{1}{4}(4a_0 - a_2^2 - a_1 a_3), \qquad b_2 = \frac{a_2}{2},$$

*and*

$$\lambda_0 = \frac{4b_2 + 2\sqrt{3b_1 + 4b_2^2}}{3}.$$

*Then the form (3) is positive semi-definite (positive definite) if and only if the $3 \times 3$ real symmetric matrix*

$$\mathbf{M}_{\lambda_0} = \begin{bmatrix} 1 & \dfrac{a_3}{2} & \dfrac{a_2 - \lambda_0}{2} \\ \dfrac{a_3}{2} & \lambda_0 & \dfrac{a_1}{2} \\ \dfrac{a_2 - \lambda_0}{2} & \dfrac{a_1}{2} & a_0 \end{bmatrix}$$

*is a positive semi-definite (positive definite) matrix.*

*Proof.* One direction of the proof is obvious since if $\lambda_0$ is a real number satisfying the matrix $\mathbf{M}_{\lambda_0}$ is positive semi-definite (positive definite), then according to the properties of positive semi-definite matrix (positive definite matrix) we have $\lambda_0 \geq a_3^2/4$ and $g(\lambda_0) = \det(\mathbf{M}_{\lambda_0}) \geq 0$ ($\lambda_0 > a_3^2/4$ and $g(\lambda_0) > 0$). Further, from Table 3, we know that some one of Case 2, Case 5, Case 6, Case 7, Case 9 (Case 2, Case 7) happens. That means that the form (3) is positive semi-definite (positive definite).

To prove the converse, we just need to show that if $\lambda_0 \geq a_3^2/4$ and $g(\lambda_0) \geq 0$, then the matrix $\mathbf{M}_{\lambda_0}$ is positive semi-definite. On one hand, if $\lambda_0 > a_3^2/4$ and $g(\lambda_0) \geq 0$, then according to Sylvester's criterion, we have $\mathbf{M}_{\lambda_0}$ is positive semi-definite. On the other hand, when $\lambda_0 = a_3^2/4$ and $g(\lambda_0) \geq 0$, since $g(a_3^2/4)$ can be represented as follows:

$$g(a_3^2/4) = \det\left(\begin{bmatrix} 1 & \dfrac{a_3}{2} & \dfrac{4a_2 - a_3^2}{8} \\ \dfrac{a_3}{2} & \dfrac{a_3^2}{4} & \dfrac{a_1}{2} \\ \dfrac{4a_2 - a_3^2}{8} & \dfrac{a_1}{2} & a_0 \end{bmatrix}\right) = -\frac{1}{256}(8a_1 - 4a_2 a_3 + a_3^3)^2,$$

we have $g(a_3^2/4) = 0$. So, we may get $a_1 = (4a_2 a_3 - a_3^3)/8$. Further, according to $\lambda_0 = (4b_2 + 2\sqrt{3b_1 + 4b_2^2})/3 = a_3^2/4$, we can also get $a_0 = (4a_2 - a_3^2)^2/64$. That means that $\mathbf{M}_{\lambda_0}$ is a positive semi-definite matrix with the rank 1.

Thus, if the form (3) is positive semi-definite, then there is no real simple intersection point between two base conics, i.e., Case 2, Case 5, Case 6, Case 7 and Case 9 in Table 1. Further, from Table 3, we have $\lambda_0 \geq a_3^2/4$ and $g(\lambda_0) \geq 0$, namely the matrix $\mathbf{M}_{\lambda_0}$ is positive semi-definite. In particular, if the form (3) is positive definite, then there is no real common point between two base conics, i.e., Case 2, Case 7. Further, from Table 3, we have $\lambda_0 > a_3^2/4$ and $g(\lambda_0) > 0$. Consequently, according to Sylvester's criterion, the matrix $\mathbf{M}_{\lambda_0}$ is positive definite. This completes the proof of the theorem. □

From the proof of Theorem 1, we can obtain directly a novel criterion for determining the positivity of a monic binary quartic form.

**Theorem 2.** *Given a monic binary quartic form*

$$f(x, y) = x^4 + a_3 x^3 y + a_2 x^2 y^2 + a_1 x y^3 + a_0 y^4. \tag{4}$$

*Let*

$$b_0 = \frac{1}{4}(-a_1^2 + a_1 a_2 a_3 - a_0 a_3^2), \qquad b_1 = \frac{1}{4}(4a_0 - a_2^2 - a_1 a_3), \qquad b_2 = \frac{a_2}{2},$$

and

$$\lambda_0 = \frac{4b_2 + 2\sqrt{3b_1 + 4b_2^2}}{3}.$$

Then the form (4) is positive semi-definite if and only if

$$\lambda_0 \geq \frac{a_3^2}{4}, \qquad -\frac{1}{4}\lambda_0^3 + b_2 \lambda_0^2 + b_1 \lambda_0 + b_0 \geq 0.$$

In particular, the form (4) is positive definite if and only if

$$\lambda_0 > \frac{a_3^2}{4}, \qquad -\frac{1}{4}\lambda_0^3 + b_2 \lambda_0^2 + b_1 \lambda_0 + b_0 > 0.$$

It is well known that the negativity of $-x^4 + a_3 x^3 y + a_2 x^2 y^2 + a_1 x y^3 + a_0 y^4$ is equivalent to the positivity of $x^4 - a_3 x^3 y - a_2 x^2 y^2 - a_1 x y^3 - a_0 y^4$. Consequently, according to the obtained conclusions, we may obtain directly the following corollaries about the negativity of a binary quartic form.

**Corollary 1.** *Given a binary quartic form*

$$f(x, y) = -x^4 + a_3 x^3 y + a_2 x^2 y^2 + a_1 x y^3 + a_0 y^4. \tag{5}$$

*Let*

$$b_1 = -\frac{1}{4}(4a_0 + a_2^2 + a_1 a_3), \qquad b_2 = -\frac{a_2}{2},$$

and

$$\lambda_0 = \frac{4b_2 + 2\sqrt{3b_1 + 4b_2^2}}{3}.$$

*Then the form (5) is negative semi-definite (negative definite) if and only if the matrix*

$$\begin{bmatrix} -1 & \frac{a_3}{2} & \frac{a_2 + \lambda_0}{2} \\ \frac{a_3}{2} & -\lambda_0 & \frac{a_1}{2} \\ \frac{a_2 + \lambda_0}{2} & \frac{a_1}{2} & a_0 \end{bmatrix}$$

*is a negative semi-definite (negative definite) matrix.*

**Corollary 2.** *Given a binary quartic form*

$$f(x, y) = -x^4 + a_3 x^3 y + a_2 x^2 y^2 + a_1 xy^3 + a_0 y^4. \tag{6}$$

*Let*

$$b_0 = \frac{1}{4}(-a_1^2 - a_1 a_2 a_3 + a_0 a_3^2), \quad b_1 = -\frac{1}{4}(4a_0 + a_2^2 + a_1 a_3), \quad b_2 = -\frac{a_2}{2},$$

*and*

$$\lambda_0 = \frac{4b_2 + 2\sqrt{3b_1 + 4b_2^2}}{3}.$$

*Then the form (6) is negative semi-definite if and only if*

$$\lambda_0 \geq \frac{a_3^2}{4}, \quad -\frac{1}{4}\lambda_0^3 + b_2 \lambda_0^2 + b_1 \lambda_0 + b_0 \geq 0.$$

*In particular, the form (6) is negative definite if and only if*

$$\lambda_0 > \frac{a_3^2}{4}, \quad -\frac{1}{4}\lambda_0^3 + b_2 \lambda_0^2 + b_1 \lambda_0 + b_0 > 0.$$

## 4. Examples

The following examples illustrate the proposed method.

*Example 1 [8].* Consider the binary quartic form

$$f(x, y) = x^4 + xy^3 + y^4.$$

According to Theorem 2, we have $b_0 = -1/4$, $b_1 = 1$, $b_2 = 0$ and further $\lambda_0 = 2/\sqrt{3}$. Since $\lambda_0 > a_3^2/4 = 0$ and $g(\lambda_0) = (16\sqrt{3} - 9)/36 > 0$, we assert that the form must be positive definite (One can check that $f(x, y) = 0$ has two pairs of complex conjugate roots in 1-dimensional projective space.).

*Example 2.* Consider the binary quartic form

$$f(x, y) = x^4 - 8x^3 y + 26x^2 y^2 - 40xy^3 + 25y^4.$$

According to Theorem 2, we have $b_0 = 1280$, $b_1 = -224$, $b_2 = 13$ and further $\lambda_0 = 56/3$. Since $\lambda_0 > a_3^2/4 = 16$ and $g(\lambda_0) = 64/27 > 0$, we assert that the form must be positive definite (One can check that $f(x, y) = 0$ has two complex conjugate double roots.).

*Example 3 [9].* Consider the binary quartic form

$$f(x, y) = x^4 + x^3 y + xy^3 + y^4.$$

According to Theorem 2, we have $b_0 = -1/2$, $b_1 = 3/4$, $b_2 = 0$ and further $\lambda_0 = 1$. Since $\lambda_0 > a_3^2/4 = 1/4$ and $g(\lambda_0) = 0$, we assert that the form is positive semi-definite (One can check that $f(x,y) = 0$ has a real double root and a pair of complex conjugate roots.).

*Example 4.* Consider the binary quartic form

$$f(x,y) = x^4 + 4x^3y + 2x^2y^2 - 4xy^3 + y^4.$$

According to Theorem 2, we have $b_0 = -16$, $b_1 = 4$, $b_2 = 1$ and further $\lambda_0 = 4$. Since $\lambda_0 = a_3^2/4 = 4$ and $g(\lambda_0) = 0$, we assert that the form is positive semi-definite (One can check that $f(x,y) = 0$ has two real double roots.).

*Example 5.* Consider the binary quartic form

$$f(x,y) = x^4 + 4x^3y + 6x^2y^2 + 4xy^3 + y^4.$$

According to Theorem 2, we have $b_0 = 16$, $b_1 = -12$, $b_2 = 3$ and further $\lambda_0 = 4$. Since $\lambda_0 = a_3^2/4 = 4$ and $g(\lambda_0) = 0$, we assert that the form is positive semi-definite (One can check that $f(x,y) = 0$ has a quadruple root.).

*Example 6.* Consider the binary quartic form

$$f(x,y) = -x^4 + 6x^3y - 13x^2y^2 + 24xy^3 - 36y^4.$$

According to Corollary 2, we have $b_0 = 0$, $b_1 = -169/4$, $b_2 = 13/2$ and further $\lambda_0 = 13$. Since $\lambda_0 > a_3^2/4 = 9$ and $g(\lambda_0) = 0$, we assert that the form is negative semi-definite (One can check that $f(x,y) = 0$ has a real double root and a pair of complex conjugate roots.).

## 5. Conclusion

In this paper, we obtained a simplified necessary and sufficient condition compared with [9] for the positivity of a binary quartic form. We considered the common points of two conics instead of roots of the given quartic form and further a novel criterion for determining the positivity was derived by using the theory of the pencil of conics. Our method required neither to use the discriminant of a quartic polynomial like the literatures [3, 5, 6, 8] nor to solve a set of inequality equations like the literature [9]. The examples illustrated the correctness of our criterion.

## Data Availability

No data were used to support the study.

## Conflict of Interest

The author declares that there is no conflict of interest regarding the publication of this paper.

# Acknowledgments

An earlier version of this paper has been presented as arXiv in Cornell University (https://arxiv.org/abs/2009.01033).